\documentclass[12pt]{article}
\usepackage{amssymb}

\def\aa{{{\scriptstyle\odot}_\nu}}
\def\A{{\mbox{\sf A}}}
\def\N{{\mbox{\sf N}}}
\def\Nn{{\mbox{\sf N}_\nu}}
\def\No{{{\mbox{\sf N}^0_\nu}}}
\def\Pol{{\mbox{\sf Pol}}}
\def\Poln{{\mbox{\sf Pol}_\nu}}
\def\SPol{{{\cal S}(\Pol)}}
\def\R{{\mathbb R}}
\def\U{{\cal U}}
\def\L{{\mbox{\sf Lin}}}
\def\sn{{\ast_\nu}}
\def\g{{\mathfrak g}}

\newtheorem{theorem}{Theorem}
\newtheorem{definition}{Definition}
\newtheorem{lemma}{Lemma}
\newtheorem{remark}{Remark}
\newtheorem{corollary}{Corollary}
\newenvironment{proof}{\noindent{\it Proof.\/}}{\vskip3mm}
\newenvironment{acknowledgement}{{\vskip3mm}
{\noindent{\bf Acknowledgements.\/}}}{\vskip3mm}

\begin{document}
\noindent{\Large\bf On Generalized Abelian Deformations}
\vskip13mm
\noindent{\sc Giuseppe  Dito\footnote{Supported by the {\em Japan
Society for the Promotion of Science} and the {\em Conseil r\'egional de
Bourgogne.}}}
\vskip2mm
\noindent{\it Research Institute for Mathematical Sciences\\ Kyoto University,
Sakyo-ku, Kyoto 606-01, Japan}
\vskip2mm
\noindent{\it D\'epartement de Math\'ematiques\\ Universit\'e de Bourgogne, 
BP 400, F-21011 Dijon Cedex, France}
\vskip5mm
\date{February 1998}

\noindent{\small \bf Abstract:} {\small We study sun-products on $\R^n$, i.e.
generalized Abelian deformations associated with star-products for general
Poisson structures on $\R^n$. We show that their cochains are given by
differential operators. As a consequence, the weak triviality of sun-products
is established and we show that strong equivalence classes are quite small.
When the Poisson structure is linear (i.e., on the dual of a Lie algebra), we
show that the differentiability of sun-products implies that covariant
star-products on the dual of any Lie algebra are equivalent each other.}

\vskip5mm 
\section{Introduction} 

A new kind of deformations was introduced in \cite{DFST97a} in connection with
the quantization of Nambu-Poisson structures (see also \cite{FDS97a}). The
main feature of these deformations is that they are not of Gerstenhaber's type
\cite{Germ64a} in the sense that one does not have a
$\mathbb{K}[[\nu]]$-algebra structure on the deformed algebra ($\mathbb{K}$ is
the ring over which is defined the original algebra $\A$ and $\nu$ denotes the
deformation parameter). More precisely, these deformations are not linear
with respect to the deformation parameter; the product operation annihilates
the deformation parameter so that one has only a $\mathbb{K}$-algebra
structure on the deformed algebra $\A[[\nu]]$.

The motivation for dealing with these generalized
deformations was that they provide non trivial Abelian deformations of the
usual product, and this point was essential for the solution proposed
in \cite{DFST97a} for the quantization of Nambu-Poisson structures. We recall
that in Gerstenhaber's framework, Abelian deformations of
the usual product of smooth functions on some manifold are always trivial (it
is a consequence of the fact that a symmetric Hochschild 2-cocycle is a
coboundary).

Explicit examples of generalized Abelian deformations were constructed in
\cite{DFST97a,DiFl97a}. There are two main classes of generalized Abelian
deformations. On the one hand, one has the Zariski products introduced in
\cite{DFST97a} which involve factorization of
polynomials in several variables into irreducible factors.
Zariski products are Abelian products on the semi-group algebra generated
by irreducible polynomials and can be constructed from any star-product
on $\R^n$. Originally the construction of a Zariski product was performed
from a Moyal product and it appeared crucial to go over semi-group algebras
with a proper notion of derivatives to fulfill algebraic requirements imposed
by the Fundamental Identity of Nambu-Poisson structures.
This construction is quite sophisticated and little is known about its
properties. Actually the Zariski quantization induced by some Zariski product
shares many properties with second quantization (appearance of a Fock space
generated by irreducible polynomials, etc.).

On the other hand, sun-products have been studied in \cite{DiFl97a}. They
have much simpler properties than Zariski products and, roughly speaking, they
can be as seen as the finite dimensional version of Zariski products. They
involve factorization into linear polynomials and can be defined on some 
algebra of functions over finite dimensional spaces.

Still generalized deformations have to find an appropriate algebraic framework
and it is the aim of this paper to study sun-products on $\R^n$ and to clarify
their structure. Our main result is that sun-products are
differentiable deformations, i.e., their cochains are differential operators
vanishing on constants. This fact allows us to find a complete
characterization of the cochains of a sun-product: Any sequence of
differential operators vanishing on linear polynomials defines a sun-product
and vice-versa.

After briefly recalling the most basic facts on star-products and Hochschild
cohomology, Sect.~2 provides a study of sun-products associated with
star-products on $\R^n$ endowed with a general Poisson structure. We show in
Theorem~\ref{diff} the differentiability of sun-products and deduce some
consequences of this property. 

We then specialize our discussion to the important case of the dual of a Lie
algebra in Sect.~3. Consider a Lie algebra $\mathfrak{g}$. Its dual
$\mathfrak{g}^*$ is endowed with a canonical Poisson structure. We show that
Gutt's star-product on $\mathfrak{g}^*$ is the only covariant star-product on
$\mathfrak{g}^*$ whose associated sun-product coincides with the usual product
on $C^\infty(\mathfrak{g}^*)$. From the differentiability of sun-products
one shows that covariant star-products on the dual of any Lie algebra are
equivalent each other.

In Sect.~4, as another consequence of the differentiable nature of
sun-products we show that sun-products are weakly trivial in the sense of
\cite{DiFl97a}. We said that two sun-products are weakly equivalent if there
exists an invertible formal series of differential intertwining these
sun-products. The sun-product operation kills all of the non-zero powers of
the deformation parameter. Weak triviality of a sun-product means weak
equivalence with the usual product (on the undeformed algebra). When one
allows the deformation parameter coming from the equivalence operator not to
be annihilated by the sun-product, one gets the notion of strong equivalence
of sun-products. By a simple argument, we remark that strong equivalence
classes are rather small. We think that these results might be helpful or give
some hints for the definition of a cohomology adapted to generalized
deformations.
\section{Sun-products on ${\mathbb R}^n$}
\subsection{Notions on star-products}
We summarize here basic facts about star-products that we shall need 
in the present paper. The general reference 
on star-products theory are the papers \cite{BFFLS78a,BFFLS78b}. 

Let $M$ be a Poisson  manifold with Poisson bracket  $P$. The space
of smooth functions $C^\infty(M)$ carries two natural algebraic structures:
It is an Abelian algebra for the pointwise product of functions and also
a Lie algebra for the Poisson bracket $P$. A star-product on $(M,P)$ is a 
formal associative deformation in the Gerstenhaber's sense \cite{Germ64a}
of the Abelian algebra structure of $C^\infty(M)$. More precisely:
\begin{definition}\label{sn}
Let $C^\infty(M)[[\nu]]$ be the space of formal series in a parameter
$\nu$ with coefficients in $C^\infty(M)$. A star-product on $(M,P)$
is a bilinear map from $C^\infty(M)\times C^\infty(M)$ to $C^\infty(M)[[\nu]]$
denoted by $f\sn g = \sum_{r\geq0}\nu^r C_r(f,g)$, $f,g\in C^\infty(M)$,
where (the cochains) $C_r\colon C^\infty(M)\times C^\infty(M)
\rightarrow C^\infty(M)$ are
bilinear maps satisfying for any $f,g,h\in C^\infty(M)$:
\begin{itemize}
\item[i)]   $C_0(f,g)=fg$; 
\item[ii)]  $C_r(c,f)= C_r(f,c)=0$, for $r\geq1$, $c\in \R$;
\item[iii)]  $\displaystyle\sum_{s+t=r\atop s,t\geq0}C_s(C_t(f,g),h) 
         \displaystyle\sum_{s+t=r\atop s,t\geq0}C_s(f,C_t(g,h))$, for
$r\geq0$;
\item[iv)] $C_1(f,g) - C_1(g,f) = 2 P(f,g)$.
\end{itemize}
\end{definition}

A star-product $\sn$ is naturally extended to a bilinear map on
$C^\infty(M)[[\nu]]$. The conditions i)--iv) above simply translate,
respectively, that a star-product is: i) a deformation of the pointwise
product; ii) it preserves the original unit ($1\sn f = f\sn 1 = f$); iii) it
is an associative product; iv) the associated star-bracket, $[f,g]_\sn =(f\sn
g - g\sn f)/2\nu$, is a Lie algebra deformation of the Lie-Poisson algebra
$(C^\infty(M),P)$.

Usually, one adds one more condition on the cochains $C_r$ of a star-product
by requiring that they should be  bidifferential operators (necessarily null
on constants by condition~ii)). These star-products are called differential
star-products. In this paper, star-product will always mean differential
star-product. One has a notion of equivalence between star-products given by:

\begin{definition}\label{equivalence}
Two star-products $\sn$ and $\sn'$ on $(M,P)$ are said to be
equivalent if there exists a formal series $T=I + \sum_{r\geq1} \nu^r T_r$,
where $I$ is the identity map on  $C^\infty(M)$ and  the $T_r$'s are
differential
operators on $C^\infty(M)$ vanishing on constants, such that 

$$ 
T(f\sn g) = T(f)\sn'T(g), \quad f,g\in C^\infty(M)[[\nu]]. 
$$ 
\end{definition} 

For a long time, star-products were known to exist on any symplectic manifold
(i.e., when the Poisson bracket $P$ is induced by some symplectic form)
\cite{DWLe83b}. Few months ago, as a consequence of his formality conjecture,
Kontsevich showed that in fact star-products exist on any Poisson manifold and
gave a complete description of their equivalence classes \cite{Konm97a}.
\subsection{Hochschild cohomology} 
Hochschild cohomology plays a  prevailing r\^ole in the deformation theory of
associative algebras. It is well known that the obstructions to equivalence of
associative deformations are in second Hochschild cohomology space and the
obstructions for extending a deformation, given up to certain order in the
deformation parameter, to the next order live in the third Hochschild
cohomology space. We shall recall here the definition and basic properties of
the Hochschild cohomology in the differentiable (null on constants) case.

Let $\A$ be the Abelian algebra $C^\infty(M)$ endowed with the pointwise
product. Consider the complex ${\cal C}^\star(\A,\A)=\{{\cal
C}^r(\A,\A)\}_{r\geq 0}$, where ${\cal C}^r(\A,\A)$ is the vector space of
$r$-linear differential operators null on constants $\phi\colon
\A^{r}\rightarrow \A$, with coboundary operator $\delta$, defined on an
$r$-cochain $C$ by:
\begin{eqnarray*}
&&\delta C (f_0,\ldots,f_r) =\cr
&&\\
&&\quad f_0 C(f_1,\ldots,f_r) + \sum_{1\leq i \leq r}
(-1)^{i}C(f_0,\ldots,f_{i-1} f_{i},\ldots,f_r) 
+ (-1)^{r+1} C (f_0,\ldots,f_{r-1})f_r,
\end{eqnarray*}
for any $f_0,\ldots,f_r$ in $\A$. The Hochschild cohomology 
(with values in $\A$) is the cohomology of the cochain complex 
$({\cal C}^\star(\A,\A),\delta)$ and shall be denoted by
${\cal H}^\star_{\rm diff,nc}(\A)$.  A fundamental result is
\begin{theorem}[Vey\cite{Veyj75a}]\label{vey}
The Hochschild cohomology ${\cal H}^\star_{\rm diff,nc}(\A)$ is isomorphic
to $\Gamma(\wedge^\star TM)$, the space of skew-symmetric contravariant
tensor fields on $M$.
\end{theorem}
Hence any Hochschild $r$-cocycle $\phi$ can be written as 
$\phi = \delta \theta + \Lambda$, where $\theta$ is an $(r-1)$-cochain 
and $\Lambda$ is an $r$-tensor on $M$. In particular, a completely 
symmetric cocycle is a coboundary.
\subsection{Notations and definitions}
We start by making precise our notations. 
The coordinates of $\R^n$ are denoted by $(x_1,\ldots,x_n)$. Let $\N$ be the
$\R$-algebra of smooth functions on $\R^n$. Let $\Pol$ be the $\R$-subalgebra
of $\N$ consisting of polynomials in $\R [x_1,\ldots,x_n]$. For a formal
parameter $\nu$, we shall denote by $\Nn$ (resp.~$\Poln$) the algebra
$\N[[\nu]]$ (resp.~$\Pol[[\nu]]$) of formal series in $\nu$ with coefficients
in
$\N$ (resp.~$\Pol$). We  distinguish in $\Nn$ a subalgebra $\No$ consisting of
formal series whose zeroth-order coefficient belongs to $\Pol$. $\Nn$, $\No$
and $\Poln$ are naturally $\R[[\nu]]$-algebras, but we shall often view them
as
$\R$-algebras.

The natural projection  $\pi\colon\Nn\rightarrow\N$ is an $\R$-algebra
homomorphism and the same symbol shall be used for the projections of $\No$
and $\Poln$ on $\Pol$.

We now define sun-products. Let $\SPol$ denote the {\em symmetric\/} tensor
algebra over $\Pol$ with symmetric tensor product $\otimes$, and let
$\lambda\colon\Pol\rightarrow\SPol$ be the $\R$-algebra homomorphism defined
by:
\begin{equation}\label{a}
\lambda(x_1^{k_1}\cdots x_n^{k_n})
=(x_1^{k_1\atop\otimes})\otimes\cdots\otimes (x_n^{k_n\atop\otimes}),
\quad \forall k_1,\ldots,k_n \geq 0.
\end{equation}
The map $\lambda$ sends a polynomial in $\Pol$ to an element of $\SPol$ by
replacing the usual product between linear factors by the symmetric tensor
product.

Let $P$ be a Poisson bracket on ${\mathbb R}^n$. Given a star-product $\sn$ on
$(\R^n,P)$, we define an $\R$-linear map  $T_\sn\colon\SPol\rightarrow\No$ by:
\begin{equation}\label{b}
T_\sn(f_1\otimes\cdots\otimes f_k)=
\frac{1}{k!}\sum_{\sigma\in S_k} f_{\sigma(1)}\sn\cdots\sn f_{\sigma(k)},
\quad \forall k\geq1,
\end{equation}
where $f_i\in\Pol$, $1\leq i\leq k$, and $S_k$ is the permutation 
group on $k$ elements.
By convention, we set $T_\sn(I)=1$, where $I$ the identity of $\SPol$.
Notice that the zeroth-order coefficient on the right-hand side of
(\ref{b}) is the product of polynomials $f_1\cdots f_k\in \Pol$, but
in general the coefficient of $\nu^r$  for $r\geq 1$ is in $\N$.
\begin{definition}\label{sunproduct}
To a star-product $\sn$ on $(\R^n,P)$, we associate a new product on
$\No$ by the following formula:
\begin{equation}\label{c}
f\aa g = T_\sn(\lambda(\pi(f))\otimes\lambda(\pi(g))),\quad
f,g\in\No. 
\end{equation} 
This product is called the $\aa$-product (or sun-product)
associated to $\sn$.
\end{definition}

In words, a sun-product on $\R^n$ associates to two polynomials $f,g\in\Pol$
the element $f\aa g \in \No$ obtained by replacing the usual product between
linear factors (in some given order) in $fg$ by a star-product $\sn$ and then
by completely symmetrizing the expression found. The extension of the product
to $f,g\in\No$ is obtained by applying the previous procedure to the
zeroth-order coefficient of $fg$. Hence a sun-product annihilates any non-zero
powers of the deformation parameter.

Basic properties of sun-products are collected in the following lemma:
\begin{lemma}
A sun-product $\aa$ on $\R^n$ is an Abelian, associative product
on $\No$. It fails to be $\R[[\nu]]$-bilinear, but it is $\R$-bilinear. 
$\No$ endowed with a product $\aa$ is an Abelian $\R$-algebra.
\end{lemma}
\begin{proof}
That the product $\aa$ is Abelian is clear from
(\ref{c}). Associativity follows from $\pi(f\aa g)= \pi(f)\pi(g)$ for
$f,g \in \No$, and from the fact that both $\lambda$ and $\pi$ are
$\R$-algebra homomorphisms: 
$f\aa(g\aa h)=T_\sn(\lambda(\pi(f))\otimes \lambda(\pi(g\aa h)))
=T_\sn(\lambda(\pi(fgh)))=(f\aa g)\aa h$, for $f,g,h\in \No$.
\end{proof}
Clearly a sun-product does not have a unit on $\No$, nevertheless one
has $1\aa f = f$ when $f$ is a linear polynomial in $\Pol$.

{}From the preceding proof, we see that to every sun-product $\aa$
we can associate a formal series of linear maps 
$\rho=\sum_{0\leq  r}\nu^r \rho_r$, where $\rho_0=Id$ is the identity
map on $\Pol$, and $\rho_r\colon \Pol\rightarrow \N$ for $r\geq1$, such that
$f\aa g = \rho(\pi(fg))$ for $f,g\in \No$. We shall (abusively) call the maps
$\rho_r$ the cochains of the sun-product $\aa$.
\subsection{Differentiability}
An example of sun-product has been explicitly computed in \cite{DiFl97a}
for  some star-product on the dual of the
Lie algebra ${\mathfrak{su}}(2)$ seen as Poisson manifold when endowed with 
its natural Lie-Poisson bracket. A remarkable feature of this sun-product is
that its cochains are differential operators. In the following, we shall show
that this fact corresponds to the general situation. As a consequence, 
any sun-product admits a natural extension from $\No$ to $\Nn$.
\begin{theorem}\label{diff}
The cochains $\rho_r$ of a sun-product $\aa$ associated
to some star-product $\sn$ on $(\R^n,P)$ are given by the restriction 
to $\Pol$ of differential operators on $\N$.
\end{theorem}
Before proving this theorem, we shall derive few lemmas. 
We consider a sun-product $\aa$  associated with some star-product $\sn$ on
$(\R^n,P)$. The cochains of the sun-product (resp. star-product) are denoted
by $\rho_r$ (resp. $C_r$). 

For any map $\phi\colon\R^k\rightarrow E$, where $E$ is a vector space, 
$\displaystyle\sum_{(i_1,\ldots,i_k)}\phi(x_{i_1},\ldots,x_{i_k})$
denotes the sum over cyclic permutations of $(x_{i_1},\ldots,x_{i_k})$.
\begin{lemma}\label{bpsi}
Let $\psi\colon \Pol\rightarrow \N$ be a linear map such that
$\psi(1)=\psi(x_i)=0$, for $1\leq i \leq n$. Let $\phi\colon
\N\times\N\rightarrow\N$ be a bidifferential operator null on constants.
If the Hochschild coboundary $\delta\psi$ satisfies for any $k\geq2$ and
indices $(i_1,\ldots,i_k)$:
\begin{equation}\label{cb}
\sum_{(i_1,\ldots,i_k)}\delta\psi(x_{i_1},x_{i_2}\cdots x_{i_k})
=\sum_{(i_1,\ldots,i_k)}\phi(x_{i_1},x_{i_2}\cdots x_{i_k}),
\end{equation}
then $\psi$ is the restriction to $\Pol$ of a differential operator 
null on constants.
\end{lemma}
\begin{proof}
On the right-hand side of Eq.~(\ref{cb}), it is clear that is sufficient
to consider bidifferential operators of the form (only these are contributing
to Eq.~(\ref{cb})):
$$
\phi(f,g)= \sum_{1\leq i \leq n} \sum_{J\atop|J|\geq1} 
\phi^{i,J} \partial_i f \partial_J g,$$
where $J=(j_1,\ldots,j_n)$ is a multi-index, $|J|=\sum_{1\leq s\leq n} j_s$,
$\partial_i = \partial/\partial x_i$, $\partial_J=\partial^{|J|}/\partial
x_1^{j_1} \cdots\partial x_n^{j_n}$ and, for fixed $i$ and $J$, $\phi^{i,J}$
is a smooth function on $\R^n$ vanishing if  $|J|$ is greater than some
integer. Consider the differential operator
$$
\tilde \psi (f) = - \sum_{1\leq i \leq n} \sum_{J\atop|J|\geq1} 
{1\over |J|+1}\phi^{i,J} \partial_{iJ} f,
$$
where $\partial_{iJ} f$ means $\partial^{|J|+1}f/\partial x_1^{j_1}
\cdots\partial x_i^{j_i +1}\cdots\partial x_n^{j_n}$ for $J=(j_1,\ldots,j_n)$.
Notice that $\tilde \psi(1)=\tilde\psi(x_i)=0$, for $1\leq i \leq n$.
The following property of $\tilde \psi$ is established by a straightforward
computation:
\begin{equation}
\sum_{(i_1,\ldots,i_k)}\delta\tilde\psi(x_{i_1},x_{i_2}\cdots x_{i_k})
=\sum_{(i_1,\ldots,i_k)}\phi(x_{i_1},x_{i_2}\cdots x_{i_k}).
\end{equation}
for any $k\geq2$ and indices $(i_1,\ldots,i_k)$.
Then, for $\psi\colon \Pol\rightarrow \N$ satisfying the hypothesis 
of the lemma, we have:
\begin{equation}\label{cpt}
\sum_{(i_1,\ldots,i_k)} 
\delta(\psi -\tilde\psi)(x_{i_1},x_{i_2}\cdots x_{i_k})=0
\end{equation}
for any $k\geq2$ and indices $(i_1,\ldots,i_k)$. 
Let $\eta = \psi -\tilde\psi$. Since
$\eta(x_{i})=0$, $1\leq i \leq n$, we have 
$\delta\eta(x_{i},f) = x_i \eta(f) - \eta(x_i f)$ for 
$1\leq i \leq n$ and $f\in \Pol$.
Then Eq.~(\ref{cpt}) implies that
$$
\eta(x_{i_1}\cdots x_{i_k})= {1\over k}\sum_{(i_1,\ldots,i_k)} 
x_{i_1}\eta(x_{i_2}\cdots x_{i_k}),
$$
and by induction on $k$, we find that $\eta=0$ on $\Pol$, 
i.e. $\psi = \tilde\psi|_{\Pol}$. $\square$
\end{proof}
\begin{lemma}\label{first}
Let $\aa$ be the sun-product associated with some 
star-product $\sn$ on $(\R^n,P)$. The first cochain $\rho_1$ of $\aa$
is a differential operator null on constants 
whose Hochschild coboundary satisfies $\delta
\rho_1=P-C_1$, where $C_1$ is the first cochain of the star-product $\sn$.
\end{lemma}
\begin{proof}
{}From Def.~\ref{sunproduct}, we have for $k\geq2$ and indices
$(i_1,\ldots,i_k)$:
\begin{eqnarray}\label{ind}
\rho(x_{i_1}\cdots x_{i_k})&=&x_{i_1}\aa\cdots\aa x_{i_k}\nonumber\\
&=&{1\over k!}\sum_{\sigma\in S_k} 
x_{i_\sigma(1)}\sn\cdots\sn x_{i_\sigma(k)},\nonumber\\
&=&{1\over k}\sum_{(i_1,\ldots,i_k)} x_{i_1}\sn\rho(x_{i_2}\cdots x_{i_k}).
\end{eqnarray}
The first-order term in $\nu$ in the last equation is:
$$
\rho_1(x_{i_1}\cdots x_{i_k}) = 
{1\over k}\sum_{(i_1,\ldots,i_k)} C_1(x_{i_1},x_{i_2}\cdots x_{i_k})
+ {1\over k}\sum_{(i_1,\ldots,i_k)} x_{i_1}\rho_1(x_{i_2}\cdots x_{i_k}),
$$
which can be written as
\begin{equation}\label{frc}
\sum_{(i_1,\ldots,i_k)} (\delta\rho_1 + C_1)(x_{i_1},x_{i_2}\cdots x_{i_k})=0,
\end{equation}
since $\rho_1(x_i)=0$. The associativity condition for a star-product implies
that $C_1$ is a Hochschild 2-cocycle and Theorem~\ref{vey} and condition iv)
in Def.~\ref{sn} tell us that $C_1 = P + \delta \theta$ where $\theta$ is a
differential operator null on constants. We can always take $\theta$ such that
$\theta(x_i)=0$, $1\leq i\leq n$, by adding a suitable 1-cocycle to it (e.g.,
$\tilde\theta(x_i) = \theta(x_i)- \sum_i \theta(x_i)\partial_i$). The Poisson
bracket $P$ is a 2-tensor and does not contribute to the left-hand side of
Eq.~(\ref{frc}). The same argument used in the proof of Lemma~\ref{bpsi} (cf.
Eq.~(\ref{cpt})) leads us to the conclusion that $\rho_1=-\theta$ and,
consequently, $\delta \rho_1=P-C_1$. $\square$
\end{proof}
{\it Proof of Theorem~\ref{diff}}. Using that the cochains of a sun-product
satisfy $\rho_r(x_i)=0$, $1\leq x_i \leq r$, we can write the equation of the
term of order $r$ in Eq.~(\ref{ind}) as:
\begin{eqnarray}\label{indu}
&&\sum_{(i_1,\ldots,i_k)} \delta\rho_r(x_{i_1},x_{i_2}\cdots
x_{i_k})=\nonumber\\
&&\nonumber\\
&&\quad -\sum_{(i_1,\ldots,i_k)} C_r(x_{i_1},x_{i_2}\cdots x_{i_k})
-\sum_{(i_1,\ldots,i_k)}\sum_{a+b=r\atop a,b\geq1}
 C_a(x_{i_1},\rho_b(x_{i_2}\cdots x_{i_k})),
\end{eqnarray}
for $k\geq2$ and $r\geq1$ (for $r=1$, the right-hand side has only one sum).
Notice that in the right-hand side of Eq.~(\ref{indu}) only the first $r-1$
cochains of the sun-product $\aa$ appear. We already know that $\rho_1$ is a
differential operator null on constants from Lemma~\ref{first}, and with the
help of Lemma~\ref{bpsi} a simple induction on $r$ proves the theorem.
$\square$
\begin{remark}
A direct consequence of Theorem~\ref{diff} is that we can extend sun-products,
originally defined on $\No$, to $\Nn$ by the formula $f\aa g = \pi(fg) +
\sum_{r\geq1}\nu^r \rho_r (\pi(fg))$ for $f,g\in \Nn$.
\end{remark}

Theorem~\ref{diff} has very simple consequences. We shall end this section by
deriving some results about the cochains of a sun-product. In Sect.~3 we shall
see that differentiability of sun-products allows one to deduce interesting
properties for star-products on the dual of a Lie algebra. The cochains of a
sun-product can be used to construct equivalence operators and this turns out
to be a quite powerful tool to establish equivalence relation between certain
type of star-products without any cohomological computations.
\begin{definition}\label{ep} 
${\cal E}(P)$ is the set of star-products on $(\R^n,P)$ such that their
associated sun-product $\aa$ coincide with the usual product on $\Pol$, i.e.
the cochains $\rho_r=0$ for $r\geq 1$. 
\end{definition} 
\begin{corollary}\label{equiep} 
Any star-product on $(\R^n,P)$ is equivalent to a star-product belonging to
${\cal E}(P)$. 
\end{corollary} 
\begin{proof} 
Let $\sn$ be a star-product and let $\{\rho_r\}_{r\geq1}$ be the cochains of
its associated sun-product. The maps $\rho_r$ are defined on $\N$ and we shall
denoted by the same symbol their $\R[[\nu]]$-linear extension to $\Nn$. Let us
define another star-product $\sn'$ by equivalence from $\sn$ with equivalence
operator $T=I +\sum_{r\geq1}\nu^r \rho_r$, that is to say: $T(f\sn' g) =
T(f)\sn T(g), \quad f,g\in \Nn.$ Since $T(x_i) = x_i$, $1\leq i\leq n$, we
have for $k\geq2$: $T(x_{i_1}\sn'\cdots\sn'  x_{i_k})= x_{i_1}\sn\cdots\sn
x_{i_k}$, and complete symmetrization  gives: $$ T(x_{i_1}\aa'\cdots\aa'
x_{i_k})= x_{i_1}\aa\cdots\aa x_{i_k}. $$ By definition $T$ is invertible and
notice that $x_{i_1}\aa\cdots\aa x_{i_k} =T(x_{i_1}\cdots x_{i_k})$, from the
equation above we conclude that $x_{i_1}\aa'\cdots\aa'  x_{i_k}= x_{i_1}\cdots
x_{i_k},$ for any $k\geq2$, i.e. the cochains of $\aa'$ satisfy $\rho_r'=0$
for $r\geq1$. Hence $\sn'$ belongs to ${\cal E}(P)$. $\square$ 
\end{proof} 
In view of the preceding corollary, the problem of classification of
equivalence classes of star-products on $(\R^n,P)$ reduces to classifying
equivalence classes in ${\cal E}(P)$. An order-by-order analysis in $\nu$ of
star-products in ${\cal E}(P)$ makes the second Lichnerowicz-Poisson
cohomology
\cite{Lica77a} space appear explicitly here. It plays the same r\^ole in the
Poisson case as the one played by the second de Rham cohomology space for the
classification of equivalences classes in the symplectic case \cite{NeTs95a}. 
\begin{corollary}
Let $\{\eta_i\}_{i\geq1}$ be a sequence of differential operators on $\N$ such
that $\eta_i(1)=\eta_i(x_k)=0$, $1\leq i,1\leq k\leq n$, and let $\sn$ be some
star-product on $(\R^n,P)$. There exists a star-product $\sn'$, equivalent to
$\sn$, such that the cochains of the sun-product $\aa'$ associated with $\sn'$
are precisely the $\eta_i$'s.
\end{corollary}
\begin{proof}
Any star-product $\sn$ is equivalent to a star-product $\sn''$ in ${\cal
E}(P)$. For $\{\eta_i\}_{i\geq1}$ satisfying the hypothesis of the corollary,
we consider a third star-product $\sn'$ defined by equivalence: $T(f\sn'' g) =
T(f)\sn' T(g)$ where $T=I +\sum_{i\geq1}\nu^i \eta_i$. It is easily verified
that the sun-product associated with $\sn'$ admits the $\eta_i$'s as cochains.
$\square$
\end{proof}
This shows that the set of possible cochains for a sun-product
on $\R^n$ coincides with the set of differential operators on $\R^n$ vanishing
on polynomial of degree less or equal to one. Also it is sufficient to
consider only one equivalence class of star-products to generate all of the
sun-products on $\R^n$. As one could have guessed, there is almost no
constraints imposed by the associativity condition on the possible cochains of
a sun-product. This fact in our opinion makes the cohomology problem for
generalized deformations quite difficult (see the discussion in Sect.~4).
\section{Sun-products on ${\mathfrak g}^*$}
We shall specialize our discussion to the case of the dual of a Lie algebra.
Let $\g$ be a real Lie algebra of dimension $n$. The dual $\g^*$ of $\g$
carries a canonical Poisson structure and, by choosing a basis of
$\g$, we can identify $\g^*$ as Poisson manifold with $\R^n$ endowed
with the following Poisson bracket:
\begin{equation}\label{poisson}
P_C(F,G)=\sum_{i,j,k=1}^n C_{ij}^k x_k 
{\partial f\over \partial x_i}{\partial g\over \partial x_j},
\quad \forall f,g\in \N,
\end{equation}
where the $C_{ij}^k$'s are the structure constants of the Lie algebra
$\g$ expressed in some basis.

A particular class of star-products which are important for physical
applications and in star-representation theory are the covariant
star-products:
\begin{definition}
Let $\g$ be a Lie algebra of dimension $n$. A star-product $\sn$
on $\R^n$ is said to be $\g$-covariant if
\begin{equation}
{1\over 2\nu}(x_i\sn x_j - x_j\sn x_i)=P_C(x_i,x_j)
=\sum_{k=1}^n C_{ij}^k x_k , \quad  \forall  1\leq i,j\leq n,
\end{equation}
where the $C_{ij}^k$'s are the structure constants of the Lie algebra $\g$ in
a
given basis.
\end{definition}
Star-products on the dual of a Lie algebra were known from the very beginning
of
the theory of star-products. The well known Moyal product is such an example,
another for ${\mathfrak{so}}(n)^*$ appears in \cite{BFFLS78b} in relation with
the quantization of angular momentum. The general case was treated by S.~Gutt
\cite{Guts83a} who defined a star-product on the cotangent bundle of any Lie
group $T^*G$. Gutt's star-product admits a restriction to $\g^*$ that we shall
call Gutt's star-product on $\g^*$.

Gutt's star-product on $\g^*$ has a simple expression that we briefly recall
here (see \cite{Guts83a} for further details). Polynomials on $\g^*$ can be
considered as elements of the symmetric algebra over $\g$, ${\cal S}(\g)$. Let
${\cal S}_r$ be the set of homogeneous polynomials of degree $r$ and let
$\U(\g)$ be the universal enveloping algebra of $\g$. The symmetrization map
$\phi\colon S(\g) \rightarrow \U(\g)$ defined by:
$$ 
\phi(X_{i_1}\cdots X_{i_r}) = \frac{1}{r!}\sum_{\sigma\in S_r}
X_{\sigma(1)}\otimes\cdots\otimes X_{\sigma(k)}, 
$$ 
(where $\otimes$ is the product in $\U(\g)$) is a bijection. Let
$\U_r=\phi({\cal S}_r)$, one has $\U(\g)=\oplus_{r\geq0} \U_r$ and each $u\in
\U(\g)$ can be decomposed as $u=\oplus_{r\geq0} u_r$, where $u_r\in \U_r$. Now
define a product between $P\in {\cal S}_p$ and  $Q\in {\cal S}_q$,  by:
$$
P\times_\nu Q = \sum_{r\geq0}(2\nu)^k
\phi^{-1}((\phi(P)\otimes\phi(Q))_{p+q-r}),
$$
and extend it by linearity to all of ${\cal S}(\g)$. It can be shown that the
product $\times_\nu$ is associative and is  defined by differential operators.
Hence one gets a star-product on ${\cal S}(\g)$ which is naturally extended to
$C^\infty(\g^*)$. This star-product is $\g$-covariant. We shall see that
Gutt's
star-product plays a special r\^ole in relation with sun-products on $\g^*$.
\begin{lemma}\label{unicity}
Let $\g$ be a fixed Lie algebra of dimension~$n$. The set of $\g$-covariant
star-products belonging to ${\cal E}(P_C)$ has only one element. In words,
there is only one $\g$-covariant star-product on $\g^*$ whose sun-product
coincides with the usual product on $\Pol$.
\end{lemma}
\begin{proof}
Let $\sn$ be a $\g$-covariant star-product on $(\R^n,P_C)$ with associated
sun-product $\aa$ which coincides with usual product on $\Pol$. Let $\L\subset
\Pol$ be the subspace of linear homogeneous polynomials on $\R^n$. It is easy
to verify that the $\sn$-powers, the $\aa$-powers and the usual powers of any
$X\in\L$ are identical:
\begin{equation}\label{power}
X^{m\atop \ast}_{} =X^{m\atop \aa}_{} = X^m_{}, \quad \forall X\in\L, m\geq0. 
\end{equation}
Obviously, we also have that $X^{m\atop \ast}_{} = X^m_{}$, for any $
X\in\L[[\nu]]$, $m\geq0$. (As usual, $\L[[\nu]]$ denotes the set of formal
series in $\nu$ with coefficients in $\L$.)

For $X\in\L[[\nu]]$, consider its $\sn$-exponential defined by:
\begin{equation}
\exp_\sn(X)=\sum_{r\geq0} {1\over r!} X^{r\atop \sn},
\end{equation}
it is an element of $\Nn$ and here $\exp_\sn(X)$ is identical to the usual
exponential $\exp(X)$ for any $X\in \L[[\nu]]$. The fact that  $\sn$ is a
$\g$-covariant star-product allows us to make usage of  the Campbell-Hausdorff
formula in the following form (in the sense of formal series):
\begin{equation}\label{CH}
\exp_\sn(sX)\sn\exp_\sn(tY)=\exp_\sn(Z(sX,tY)), \quad X,Y \in \L, s,t\in \R,
\end{equation}
where $Z(X,Y)=\sum_{r\geq0} \nu^r Z_r(X,Y) \in\L[[\nu]]$, and the $Z_r$'s are
related to the Campbell-Hausdorff coefficients by $Z_r(X,Y)= 2^r c_{r+1}(X,Y)$
(where $c_1(X,Y) = X+Y$, $c_2(X,Y) = P_C(X,Y)/2$, etc.). As the
$\sn$-exponential of $X\in\L[[\nu]]$ is simply the usual exponential,
Eq.~(\ref{CH}) yields
\begin{equation}\label{CHs}
\exp(sX)\sn\exp(tY)=\exp(Z(sX,tY)), \quad X,Y \in \L,s,t\in \R,.
\end{equation}
Hence a $\g$-covariant star-product for which the associated $\aa$-product is
the usual product must satisfy the preceding relation.

Actually Eq.~(\ref{CHs}) determines the star-product $\sn$ completely. Notice
that a bidifferential operator $B\colon \N\times\N\rightarrow \N$ is
completely
characterized by the functions $B(X^a,Y^b)$, $a,b\in {\mathbb N}, X,Y \in \L$.
The functions $C_r(X^a,Y^b)$, $0\leq a,b,r$, $X,Y \in \L$, which completely
determine the cochains $C_r$ of $\sn$ can be easily computed by
differentiation
with respect to $s$ and $t$ on both sides of Eq.~(\ref{CHs}) of the
coefficient
of $\nu^r$ and by evaluation at $s=t=0$. Therefore there is at most one
star-product whose associated $\aa$-product is the usual product on $\Pol$.

It is easy to show that the star-product defined by Eq.~(\ref{CHs}) has a the
usual product as associated $\aa$-product. By setting  $Y=X$ in
Eq.~(\ref{CHs}), we find that $X^a\sn X^b =X^{a+b}$, $\forall X\in\L$, $0\leq
a,b$, which implies by induction that $X^{a\atop\sn} = X^a$, $\forall X,a,b$.
By Eq.~(\ref{power}) we have  $X^{a\atop\aa} = X^a$, $\forall X,a$, and since
$\aa$ is Abelian, it implies that $\aa$ is the usual product on $\Pol$ and
this
proves the lemma. $\square$
\end{proof}
\begin{lemma}\label{gutt}
Let $\g$ be a Lie algebra. The $\g$-covariant star-product characterized in
Lemma~\ref{unicity} is Gutt's star-product on $\g^*$.
\end{lemma}
\begin{proof}
We shall use the notations introduced in the proof of Lemma~\ref{unicity}. Let
$\g$ be of dimension $n$ and let $\sn$ be the star-product characterized in
Lemma~\ref{unicity} by Eq.~(\ref{CHs}). The identification of the coefficients
of $\nu^r$ in (\ref{CHs}) gives:
\begin{equation}\label{cochain}
C_r(\exp(sX),\exp(tY)) = F_r(sX,tY)\exp(sX+tY),
\quad \forall X,Y\in \L, s,t\in \R,
\end{equation}
where the $F_r$'s are polynomial functions of the (normalized)
Campbell-Hausdorff coefficients $Z_r(sX,tY)$ and are defined by the following
recursive relation with $F_0=1$:
\begin{equation}\label{fr}
F_r = {1\over r}\sum_{k=0}^{k=r-1}(r-k)Z_{r-k} F_k, \quad r\geq1.
\end{equation}
By induction, one finds the explicit expression for $F_r$ for $r\geq1$ to be:
\begin{equation}\label{efr}
F_r =\sum_{k=1}^{k=r}
\sum_{m_1>\cdots>m_k\geq1\atop{ n_1,\ldots n_k \geq 1\atop
m_1 n_1+\cdots+ m_k n_k=r}}
{1\over n_1!\cdots n_k !}(Z_{m_1})^{n_1}\cdots (Z_{m_k})^{n_k}.
\end{equation}
Now we shall derive an explicit expression for $X\sn \exp(Y)$, $X,Y\in\L$. 
Notice that this relation also characterizes $\sn$ as any polynomial can be
expressed as a $\sn$-polynomial 
(it is a simple consequence of Eq.~(\ref{power})).
In general, the Campbell-Hausdorff coefficients 
$\{c_i\}_{i\geq1}$ ($c_1(X,Y) = X+Y$, 
$c_2(X,Y) = {1\over2}[X,Y]$, etc.) have the following properties:
\begin{eqnarray}\label{chpro}
&&c_i(0,X) = c_i(X,0) =0,\quad i\geq 2;\nonumber \\
&&\nonumber \\
&&{\partial\over\partial s}c_i(sX,Y)|_{s=0}
= {B_{i-1}\over (i-1)!}(ad_Y)^{i-1}(X),\quad i\geq 2;
\end{eqnarray}
where $ad_Y\colon X\mapsto [Y,X]$, and $B_n$ are the Bernoulli numbers.
These can be easily derived from the standard recursive formula for the
$c_i$'s,
see e.g. \cite{NaSt82a}.

Also, using Eqs.~(\ref{efr}) and (\ref{chpro}), 
along with the definition of $Z_r$ ($Z_r = 2^r c_{r+1}$),
one finds that
\begin{eqnarray*}
&&F_r(0,Y)=0, \quad r\geq1;\\
&&\\
&&{\partial\over\partial s}F_r(sX,Y)|_{s=0}
= {\partial\over\partial s}Z_r(sX,Y)|_{s=0}=
{2^rB_{r}\over r!}(ad_Y)^{r}(X),\quad r\geq 1.
\end{eqnarray*}
Therefore we can write
\begin{eqnarray}\label{eco}
C_r(X,\exp(tY)) 
&=& {\partial\over\partial s}(F_r(sX,tY)\exp(sX+tY))|_{s=0}\nonumber\\
&&\nonumber\\
&=& {2^r{B_{r}}\over r!}(ad_{tY})^{r}(X)\exp(tY),\quad r\geq 1.
\end{eqnarray}
For $r=0$, we simply have: $C_r(X,\exp(tY)) = X\exp(tY)$. Equation~(\ref{eco})
is also characterizing Gutt's star-product on $\g^*$
(compare with Eq.~(3.2) in \cite{Guts83a}). $\square$
\end{proof}
As a simple consequence of Lemmas~\ref{unicity} and \ref{gutt} we
have the following corollary which  tells us that any two covariant
star-products on the dual of a Lie algebra are equivalent.
\begin{corollary}\label{coequi}
Any covariant star-product on the dual of a Lie algebra $\g$
is equivalent to Gutt's star-product on $\g^*$.
\end{corollary}
\begin{proof}
Let $\sn$ be a $\g$-covariant star-product on $\g^*$, the dual of a Lie
algebra of dimension~$n$. From Corollary~\ref{equiep}, $\sn$ is equivalent to
a star-product $\sn'$ belonging to ${\cal E}(P_C)$, where $P_C$ is the
Lie-Poisson structure on $\g^*$. The equivalence operator is constructed out
from the cochains of the sun-product associated with $\sn$ and it leaves
invariant linear polynomials, i.e. $T(x_i)=x_i$, $1\leq i\leq n$.
Consequently, $\sn'$ is also a $\g$-covariant star-product. According to
Lemmas~\ref{unicity} and~\ref{gutt}, $\sn'$ must be Gutt's star-product on
$\g^*$. $\square$
\end{proof}
\begin{remark}
Though the de Rham cohomology of $\g^*$ is trivial, not all star-products on
$\g^*$ are equivalent. Indeed, in the symplectic case, the second de Rham
cohomology space classifies equivalence classes of star-products. In the
Poisson case, one has to consider the Lichnerowicz-Poisson cohomology
\cite{Lica77a} instead, and this cohomology is not in general trivial for the
Lie-Poisson structure on $\g^*$. See \cite{Meld89a}, for explicit computations
of some of the (Chevalley-Eilenberg) cohomology spaces for the dual of a Lie
algebra.
\end{remark}
\section{Weak and strong equivalences}
In the deformation theory of some algebraic structure one has the notion of
equivalent deformations. The equivalence of star-products given
Def.~\ref{equivalence} is adapted to the associative (differential) case and
one has similar notions of equivalence for other algebraic structures (e.g.,
Lie algebras, Abelian algebras, etc.). Moreover, as mentioned in Sect.~2.2, it
is a general result of Gerstenhaber \cite{Germ64a} that obstructions for
equivalence of deformations reside in the second cohomology space of an
appropriate cohomology. For associative, Lie, Abelian deformations the
associated cohomologies are, respectively, Hochschild, Chevalley-Eilenberg,
Harrison cohomologies. One may wonder what is the corresponding  cohomology
for generalized Abelian deformations. Before discussing on that matter, it is
important to bear in mind that in Gerstenhaber's theory of deformations a
deformed algebraic structure has a structure of ${\mathbb K}[[\nu]]$-algebra,
where ${\mathbb K}$ is the ground ring of the original structure. This
feature, which is crucial to determine the appropriate cohomology, does not
hold anymore in the case of generalized deformations.

The answer to the cohomology issue raised by generalized deformations might
be, as advocated by M.~Flato \cite{Flato96}, that one has to give a
noncommutative ring structure on the space of formal parameters in such a way
that $\R[[\nu]]$-bilinearity would be restored. This should lead to a
noncommutative deformation theory and the first steps toward this program were
taken by Pinczon \cite{Ping97b} who considered the case where the deformation
parameter is acting by different left and right endomorphisms on the algebra
(hence the deformation parameter is not required to commute with the
undeformed algebra). This point of view produced very interesting results
(e.g., deformation of the Weyl algebra yields supersymmetric algebras), but
still generalized deformations do not fit in the particular framework
considered in \cite{Ping97b}.

The cohomology problem is still open and in a previous work \cite{DiFl97a}
we have nevertheless considered two notions of
equivalence for sun-products. They are mimicking the usual notion of
equivalence and take into account that sun-products are not
$\R[[\nu]]$-bilinear operations, but only $\R$-bilinear. Let us recall their
definitions.

\begin{definition}\label{equisn}
Two sun-products $\aa$ and $\aa'$ on $(\R^n,P)$ are said to be (a) weakly
((b) strongly) equivalent, if there exists an
${\mathbb R}[[\nu]]$-linear map $S_\nu\colon
\Nn\mapsto \Nn$ where $S_\nu=\sum_{r\geq 0} \nu^r S_r$,
with $S_r\colon \N\rightarrow\N$, $r\geq1$, being differential operators
and $S_0=I$, such that for $f,g \in \N$ the following holds:
\begin{eqnarray*}
(a)\ &S_\nu(f \aa g) = &S_\nu(f) \aa' S_\nu(g),\\
&&\\
(b)\ &S_\nu(f \aa g) = &S_\mu(f) \aa'S_\mu(g)|_{\mu=\nu}.
\end{eqnarray*}
\end{definition}
For weak equivalence, condition  (a) above can be equivalently replaced by
$S_\nu(f \aa g) = f \aa' g$, as sun-products annihilate the deformation
parameter $\nu$. In the case of strong equivalence, condition (b), when
written in terms of the cochains of the sun-products, simply states that:
\begin{equation}\label{strong}
\sum_{r+s=t\atop r,s\geq0}S_r(\rho_s(fg))=
\sum_{r+a+b=t\atop r,a,b\geq0}\rho_r'(S_a(f) S_{b}(g)),\quad f,g\in\N,t\geq0,
\end{equation}
where the $\rho_i$'s (rep.~$\rho_i'$'s) are the cochains of $\aa$
(resp.~$\aa'$). It can be easily checked that Def.~\ref{equisn} indeed defines
equivalence relations on the set of  sun-products. Weak or strong triviality
has to be understood as weak or strong equivalence with the pointwise  product
on $\N$.

We shall now draw some conclusions for weak and strong equivalences of
sun-products from Theorem~\ref{diff}. It was shown in \cite{DiFl97a} that a
sun-product is weakly trivial if its cochains are differential operators.
Hence as a corollary of Theorem~\ref{diff}, we simply have:
\begin{corollary}
Let $\aa$ be a  sun-product on $(\R^n,P)$, then $\aa$ is weakly trivial.
\end{corollary}
\begin{proof} Let $\rho_i$ be the cochains of $\aa$. They are
differential operators null on constants by Theorem~\ref{diff}. Then define
$S_\nu$ to be the formal inverse of $\sum_{r\geq0}\nu^r \rho^r$. The map
$S_\nu$ satisfies $S_\nu(f\aa g)=f\cdot g$ for $f,g\in \N$, where $\cdot$
denotes the pointwise product, hence $\aa$ is weakly equivalent to the
pointwise product. $\square$ \end{proof}
On the the hand, we shall see that strong equivalence puts severe conditions
on the equivalence operator $S_\nu$. By setting $g=1$ in Eq.~(\ref{strong}),
we get with shortened notations that $S_\nu \rho =\rho' S_\nu$ and by
substituting this relation in Eq.~(\ref{strong}), we find that the equivalence
operator should satisfy $S_\nu(fg)=S_\nu(f)S_\nu(g)$. Hence $S_\nu$ can
be nothing else than the exponential of a formal series of derivations of the
pointwise product. Actually there are still some supplementary constraints on
$S_\nu$, but we do not need to be concerned with them. We conclude that strong
equivalence classes are very small and can even reduce to a single point
in some situations (e.g., the equivalence class of strongly trivial
sun-products).

Although, we do not know whether weak and strong equivalences are induced by
the cohomology of some complexes, these notions provide limiting cases
between which a proper notion of equivalence for generalized deformations
should lie.
\begin{acknowledgement}
The author would like to thank Mosh\'e Flato and Daniel Sternheimer for very
useful discussions, and Izumi Ojima for great hospitality at RIMS where this
work was finalized.
\end{acknowledgement}

\end{document}